\documentclass[reqno,12pt]{amsart}
\usepackage{mathrsfs}
\usepackage{txfonts}
\usepackage{amssymb,latexsym}
\usepackage{amsmath,amsfonts,amsthm}
\usepackage{graphicx}
\usepackage{indentfirst}
\usepackage{fancyhdr}
\usepackage{CJK}
\usepackage[margin=3cm]{geometry}
\usepackage[T1]{fontenc}
\allowdisplaybreaks[4]

\fancypagestyle{plain}{
 \fancyhead[LO,RE]{\usebox{\mygraphic}}
 
 }
\begin{document}
\setlength{\abovedisplayshortskip}{3mm}
\setlength{\belowdisplayshortskip}{3mm}
\setlength{\abovedisplayskip}{3mm}
\setlength{\belowdisplayskip}{3mm}

\newtheorem{theorem}{Theorem}[section]
\newtheorem{corollary}[theorem]{Corollary}
\newtheorem{lemma}[theorem]{Lemma}
\newtheorem{definition}[theorem]{Definition}
\newtheorem{proposition}[theorem]{Proposition}
\newtheorem{remark}[theorem]{Remark}
\newtheorem{example}[theorem]{Example}
\newtheorem{theoremalph}{Theorem}
 \renewcommand\thetheoremalph{\Alph{theoremalph}}
\newcommand{\bta}{\begin{theoremalph}}
\newcommand{\eeta}{\end{theoremalph}}
\newcommand{\bth}{\begin{theorem}}
\newcommand{\eeth}{\end{theorem}}
\newcommand{\ble}{\begin{lemma}}
\newcommand{\ele}{\end{lemma}}
\newcommand{\bco}{\begin{corollary}}
\newcommand{\eco}{\end{corollary}}
\newcommand{\bde}{\begin{definition}}
\newcommand{\ede}{\end{definition}}
\newcommand{\bpr}{\begin{proposition}}
\newcommand{\epr}{\end{proposition}}
\newcommand{\bre}{\begin{remark}}
\newcommand{\ere}{\end{remark}}
\newcommand{\beg}{\begin{example}}
\newcommand{\eeg}{\end{example}}

\newcommand{\beq}{\begin{equation}}
\newcommand{\eeq}{\end{equation}}
\newcommand{\ben}{\begin{equation*}}
\newcommand{\een}{\end{equation*}}
\newcommand{\beqn}{\begin{eqnarray}}
\newcommand{\eeqn}{\end{eqnarray}}
\newcommand{\be}{\begin{eqnarray*}}
\newcommand{\ee}{\end{eqnarray*}}
\newcommand{\ban}{\begin{align*}}
\newcommand{\ean}{\end{align*}}
\newcommand{\bal}{\begin{align}}
\newcommand{\eal}{\end{align}}
\newcommand{\bad}{\aligned}
\newcommand{\ead}{\endaligned}
\newcommand{\lan}{\langle}
\newcommand{\ran}{\rangle}

\newcommand{\na}{\nabla}
\newcommand{\vp}{\varphi}
\newcommand{\La}{\Lambda}
\newcommand{\la}{\lambda}
\newcommand{\Om}{\Omega}
\newcommand{\ta}{\theta}
\newcommand{\fr}{\frac}
\newcommand{\iy}{\infty}
\newcommand{\ve}{\varepsilon}
\newcommand{\pa}{\partial}
\newcommand{\al}{\alpha}
\newcommand{\mr}{\mathbb{R}^n}
\newcommand{\bu}{\bullet}
\newcommand{\si}{\sigma}

\newenvironment{sequation}{\begin{equation}\small}{\end{equation}}
\newenvironment{tequation}{\begin{equation}\tiny}{\end{equation}}

\title
{The concavity of $p$-entropy power and applications in functional inequalities}
\author
{Yu-Zhao Wang}
\author{Xinxin Zhang}
\address{School of Mathematical Sciences, Shanxi University, Taiyuan, 030006, Shanxi, China}
\email{wangyuzhao@sxu.edu.cn}
\email{zhangxinxin.mail@qq.com}

\maketitle
\maketitle \numberwithin{equation}{section}
\maketitle \numberwithin{theorem}{section}
\setcounter{tocdepth}{2}
\setcounter{secnumdepth}{2}

\begin{abstract}
In this paper, we prove the concavity of $p$-entropy power of probability densities solving the $p$-heat equation on closed Riemannian manifold with nonnegative Ricci curvature. As applications, we give new proofs of $L^p$-Euclidean Nash inequality and $L^p$-Euclidean Logarithmic Sobolev inequality, moreover, an improvement of $L^p$-Logarithmic Sobolev inequality is derived.

\vspace{5mm}
\textbf{Mathematics Subject Classification (2010)}. Primary 35K92; Secondary 58J35.
\vspace{5mm}

\textbf{Keywords}. Entropy power inequality, Concavity, $p$-heat equation, Logarithmic Sobolev inequality,  Nash inequality.

\end{abstract}


\section{Introduction and main results}

In 1948 classic paper \cite{Shannon}, Shannon introduced the entropy power,
\begin{equation*}
N(X)=\exp\left(\frac2nH(X)\right).
\end{equation*}
and gave a proof of entropy power inequality(EPI) \eqref{EPI} for independent random variables $X$ and $Y$,
\begin{equation}\label{EPI}
N(X+Y)\ge N(X)+N(Y),
\end{equation}
where $H(X)$ was the entropy for random variable $X$ with a probability density function $u$,
\begin{equation*}
H(X)=-\int_{\mathbb{R}^{n}}u\log u\,dx.
\end{equation*}
The entropy power inequality plays an important role in the fields of information
theory, probability theory and convex geometry. There are many interesting consequences on the entropy power inequality, one of these consequences has been noticed by Costa. More precisely, in 1985, Costa \cite{Costa} proved that, if $u(t),\,t > 0$, are probability densities solving the heat equation $\partial_tu=\Delta u$ in the whole space $\mathbb{R}^n$, then
\begin{equation}\label{CEP}
\frac{d^2}{dt^2}{N}(u(t))\le0.
\end{equation}
Inequality \eqref{CEP} is referred to as the \textbf{concavity of entropy
power}. In \cite{Villani}, Villani gave a direct proof of \eqref{CEP} in a strengthened version with an exact error term, which connects the concavity of entropy power with some identities of Barky-\'Emery.

Recently, G.Savar\'e and G.Toscani \cite{ST} showed
that the concavity of entropy power
is a property which is not restricted to Shannon entropy power
 in connection with the heat equation, but it holds for
the $\gamma$-th R\'enyi entropy power \eqref{RenPow},
\begin{equation}\label{RenPow}
{N}_{\gamma}(u)\doteqdot\exp\left(\frac{\lambda}{n}{R}_{\gamma}(u)\right),\quad \lambda=2+n(\gamma-1)>0,
\end{equation}
which connects
with the solution to the nonlinear diffusion equation
\begin{equation}\label{NDE}
\partial_tu=\Delta u^{\gamma},
\end{equation}
where ${R}_{\gamma}$ is the R\'enyi entropy
\begin{equation*}
{R}_{\gamma}(u)\doteqdot\frac1{1-\gamma}\log\int_{\mathbb{R}^n}u^{\gamma}(x)\,dx,\quad\gamma\in(0,\infty), \;\gamma\neq1.
\end{equation*}
When $\gamma > 1-\frac2n$, they have proved the concavity of R\'enyi
entropy power defined in \eqref{RenPow},
\begin{equation*}
\frac{d^2}{dt^2}{N}_{\gamma}(u(t))\le0,
\end{equation*}
where $u(t),t>0$ are probability densities solving \eqref{NDE}in $\mathbb{R}^n$.

In this paper, motivated by above works, we study the entropy with respect to nonlinear diffusion on $\mathbb{R}^n$ and Riemannian manifolds. Let $u$ be a positive solution to the $p$-heat equaiton
\begin{equation}\label{pheat}
\frac{\partial u^{p-1}}{\partial t}=(p-1)^{p-1}\Delta_{p}u,
\end{equation}
where $\Delta_{p}u={\rm div}(|\nabla u|^{p-2}\nabla u)$ is the $p$-Laplacian of $u$,
define the $p$-entropy
\begin{equation}\label{pentropy}
H_p(u)\doteqdot-\int_Mu^{p-1}\log u^{p-1}\,dV,\quad p>1
\end{equation}
and $p$-entropy power
\begin{equation}\label{pSEP}
N_p(u)\doteqdot\exp\left(\frac pn{H}_p(u)\right)
\end{equation}
on Riemannian manifold $M$, so that the Shannon's entropy and entropy power can be identified with the $p$-entropy and $p$-entropy power of index $p=2$ respectively.

Kotschwar and Ni \cite{Ni2009} introduced $p$-entropy \eqref{pentropy} on Riemannian manifold and proved the Perelman type $W$-entropy monotonicity formula with nonnegative Ricci curvature. The first author generalized this result to the weighted Riemannian manifold with nonnegative $m$-Bakry-Emery Ricci curvarure \cite{WYC} and $m$-Bakry-Emery Ricci curvarure bounded below \cite{WYZ} respectively.

The first result in this paper is the concavity of $p$-entropy power with respect to the $p$-heat equation on closed Riemannian manifold with nonnegative Ricci curvature. Due to this property and motivated by Villani \cite{Villani}, we establish a deep link between entropy, information and entropy power by means of an nonnegative quantity. We introduce the $p$-Fisher information
\begin{equation}\label{pfisher}
\frac{d}{dt}H_{p}(u)=I_{p}(u)=(p-1)^{p}\int_{M}\frac{|\nabla u|^{p}}{u}\,dV
\end{equation}
for the second order derivative of the $p$-entropy power resorting to the DeBruijn's identity.
The precise argument is following.

\begin{theorem}\label{Error}
Let $p>1$ and $u(x,t),t>0$ be a positive solution to the $p$-heat equation \eqref{pheat} on closed Riemannian manifold $(M,g)$. Assume $v=-(p-1)\log u$ and $\omega=|\nabla v|^{2}$, then we have
\begin{equation}
\frac{d^{2}}{dt^{2}}N_{p}(u)=-\frac{p^2}{n}N_{p}(u)\int_M\left(\left|\omega^{\frac{p}{2}-1}\nabla_{i}\nabla_{j}v
-\frac{I_{p}(u)}{n}a_{ij}\right|_{A}^{2}+\omega^{p-2}{\rm Ric}(\nabla v,\nabla v)\right) e^{-v}\,dV,
\end{equation}
where $I_p(u)$ is the $p$-Fisher information defined in \eqref{pfisher} with respect to $p$-entropy and for any $2$-tensor $T$, $|T|^2_A=a^{ij}a^{kl}T_{ik}T_{jl}$, $a^{ij}=g^{ij}+(p-2)\frac{v^{i}v^{j}}{\omega}$ is the inverse of $(a_{ij})$. When the Ricci curvature is nonnegative, the $p$-entropy power defined in \eqref{pSEP} is concave.
\end{theorem}
\begin{remark} When $M$ is an Euclidean space $\mathbb{R}^n$ and $u(x,t)$ is a smooth and rapidly decaying positive solution to the $p$-heat equation, then
\begin{equation}\label{concave}
\frac{d^{2}}{dt^{2}}N_{p}(u)=-\frac{p^2}{n}N_{p}(u)\int_{\mathbb{R}^n}\left|\omega^{\frac{p}{2}-1}\nabla_{i}\nabla_{j}v
-\frac{I_{p}(u)}{n}a_{ij}\right|_{A}^{2}e^{-v}\,dx\le0,
\end{equation}
that is the $p$-entropy power defined in \eqref{pSEP} is concave.
\end{remark}
 In a recent preprint \cite{LL}, S.Z. Li and X-.D. Li studied the entropy power on the weighted Reimannian manifold, they proved that the Shannon entropy power and R\'enyi entropy power are concave under the nonnegative Barky-\'Emery Ricci curvature condition \cite{Bakry}. Motivated by their work, we can also prove the concavity of $p$-entropy power with nonnegative Barky-\'Emery Ricci curvature in another paper \cite{WZ2}.

There is an explicit connection between the $p$-entropy power and the solution to $p$-heat equation, which can be highlighted owing to its fundamental solution. In the preceding result, we assume $u(x,t)$ as a solution to $p$-heat equation \eqref{pheat} and obtain the fact that the second order derivative of $p$-entropy power is non-positive. Indeed, the case equal to zero can be achieved by the fundamental solution to \eqref{pheat}, which takes the form
 \begin{equation}\label{Gpt}
 G_{p,t}(x)=t^{-\frac{n}{p(p-1)}}\widetilde{G}_{p}\left(\frac{x}{t^{\frac{1}{p}}}\right)
 \end{equation}
from the time-independent function
\begin{equation}\label{Gp}
\widetilde{G}_{p}(x)=\left(C_{p,n}e^{-\varphi_{0}(x,1)}\right)^{\frac{1}{p-1}}
\end{equation}
with
$$C_{p,n}\doteqdot(p^{\frac{1}{p}}q^{\frac{1}{q}})^{-n}\pi^{-\frac{n}{2}}
\frac{\Gamma(\frac{n}{2}+1)}{\Gamma(\frac{n}{q}+1)},
\;\varphi_{0}(x,t)\doteqdot\frac{p-1}{p^{q}}\left(\frac{| x|}{t^{\frac{1}{p}}}\right)^{q},\;\frac{1}{p}+\frac{1}{q}=1,\;x\in \mathbb{R}^{n},\;t>0.
$$
As for $p>1$, direct calculations indicate that $p$-entropy power of the fundamental solution \eqref{Gpt} is linear on time $t$, that is,
\begin{equation}\label{tGpt}
N_{p}(G_{p,t}(x))=t\cdot N_{p}(\widetilde{G}_{p}(x))
\end{equation}
with the result $\frac{d^{2}}{dt^{2}}N_{p}(G_{p,t}(x))=0$.
This expression of formula \eqref{Gpt}, which involves a couple of Gaussian densities for the choice of $p=2$ with variance equal to $2$, $2t$ respectively, is close to the definition of self-similar solution to the nonlinear diffusion equation \eqref{NDE} (\cite{ST}), in which we observe the same property as equation \eqref{tGpt}, and so that highlight the point that $p$-entropy power is rigorously concave apart from its fundamental solution \eqref{Gpt} of equation \eqref{pheat}.

The concavity ensure that a functional (see \eqref{FI}) with respect to the first derivative of $p$-entropy power is monotonously non-increasing. Thus it will reach its lower bound along the solutions to $p$-heat equation \eqref{pheat} as time tends to infinity. In this note, let us introduce the second result below.

\begin{theorem}\label{isoper}
Let $p>1$, $\| u^{p-1}\|_{L^{1}}=1$, every smooth and rapidly decaying positive function $u(x)$ satisfies
\begin{equation}\label{isoperi}
N_{p}(u)I_{p}(u)\geq
N_{p}\left(\widetilde{G}_{p}(x)\right)I_{p}\left(\widetilde{G}_{p}(x)\right)=\gamma_{n,p},
\end{equation}
where the value of the strictly positive constant $\gamma_{n,p}$ is given by
\begin{equation}\label{gammacon}
\gamma_{n,p}=n(qe)^{p-1}\pi^{\frac{p}{2}}\left[\frac{\Gamma(\frac{n}{2}+1)}{\Gamma(\frac{n}{q}+1)}\right]^{-\frac{p}{n}},\quad q=\frac p{p-1}.
\end{equation}
\end{theorem}
Theorem \ref{isoper} is a well-known property which goes back to the concavity of $p$-entropy power in the sharp form \eqref{isoperi}. When $p=2$, inequality \eqref{isoperi} turns to be
\begin{equation}\label{IEPI}
N_{2}(u)I_{2}(u)\geq2\pi ne,
\end{equation}
which known as {\bf Isoperimetric entropy power inequality}. In \cite{Toscani}, G.Toscani proved the logarithmic Sobolev inequality and the Nash's inequality
with the sharp constant again by using of isoperimetric entropy power inequality \eqref{IEPI}.
In \cite{ST}, G.Savar\'e and G.Toscani implied an improvement of Sobolev inequality by means of rewriting the isoperimetric inequality for entropy. Later, the concavity of entropy power and isoperimetric entropy power inequality can be viewed as  powerful tools to prove various sharp functional inequalities. Motivated by their works, we can also show that some inequalities as applications of concavity theorem \eqref{concave} and isoperimetric inequality \eqref{isoperi}. The first application is the $L^p$-Nash inequality.
\begin{theorem}\label{OPNash} For $p>1$ and $q=\frac p{p-1}$, assume $g(x)$ is a positive function in $L^p(\mathbb{R}^n)$, then
\begin{equation}\label{PNashineq}
\left(\int_{\mathbb{R}^{n}}g^{p}dx\right)^{1+\frac{q}{n}}\leq
p^{p}\gamma^{-1}_{n,p}\left(\int_{\mathbb{R}^{n}}gdx\right)^{\frac{pq}{n}}\int_{\mathbb{R}^{n}}|\nabla
g|^{p}\,dx
\end{equation}
where $\gamma_{n,p}$ is the constant in \eqref{gammacon}.
\end{theorem}
\begin{remark}
When $p=2$, the inequality reduces \eqref{PNashineq} to the classical Nash's inequality in sharp form
\begin{equation}
\left(\int_{\mathbb{R}^{n}}g^{2}\,dx\right)^{1+\frac{2}{n}}\le\frac{2}{\pi en}\left(\int_{\mathbb{R}^{n}}g\,dx\right)^{\frac{4}{n}}\int_{\mathbb{R}^n}|\nabla g|^2\,dx.
\end{equation}
\end{remark}

The second application is a new proof of $L^p$-Logarithmic-Sobolev inequality by the concavity of $p$-entropy power. Furthermore, we give an improvement of $L^p$-Logarithmic-Sobolev inequality based on a generalized Csisz\'ar-Kullback inequality.
\begin{theorem}\label{improve}
Let $p>1$, $\theta\in(0,2]$ and $g(x)$ be a positive function. Then we have
\begin{equation}\label{LpEu}
\int_{\mathbb{R}^{n}}g^{p}\log g^{p}\,dx\leq\frac{n}{p}\log\left(p^{p}\gamma^{-1}_{n,p}\int_{\mathbb{R}^{n}}|\nabla g|^{p}\,dx\right)
\end{equation}
and
\begin{equation}\label{GLogSobolev}
p^{p}\int_{\mathbb{R}^{n}}|\nabla g|^{p}\,dx-\left(\int_{\mathbb{R}^{n}}g^{p}\log g^{p}\,dx
+\frac{n}{p}\log\left(\frac{pe}{n}\gamma_{n,p}\right)\right)
\geq
\frac{p}{8n}\left\{\int_{\mathbb{R}^{n}}\left(g^{p}/\widetilde{G}_{p}^{p-1}-1\right)^{\theta}\widetilde{G}_{p}^{p-1}\,dx\right\}^{\frac{4}{\theta}},
\end{equation}
where $\gamma_{n,p}$ is a constant in \eqref{gammacon} and $\widetilde{G}_{p}$ is defined in \eqref{Gp}.
\end{theorem}

\vspace{3mm}
An outline of this paper is as follows. Based on $p$-heat equation \eqref{pheat}, we show the concavity of $p$-entropy power(Section 2) and the corresponding isoperimetric entropy power inequality (section 3).  In section 4 and 5, we obtain some applications of $p$-entropy power, which are generalizations of the $L^p$-Euclidean Nash inequality and (improved) $L^p$-Euclidean logarithmic Sobolev inequality.

\section{The concavity of $p$-entropy power}

Before starting our proof, we need following two dissipation formulae.

\begin{lemma}\label{lemma}
Let $p>1$ and $u(x,t),t>0$ be a positive solution to the
$p$-heat equation \eqref{pheat} on Riemannian manifold $(M,g)$. Assume $v=-(p-1)\log u$, and set $\omega=|\nabla v|^{2}$, then $p$-entropy defined in \eqref{pentropy}
can be rewritten in its simplified form
\begin{equation}\label{pentropy1}
H_{p}(u)=\int_{M}ve^{-v}\,dV.
\end{equation}
Hence,
\begin{equation}
\frac{d}{dt}H_p(u)=\int_{M}\omega^{\frac{p}{2}}e^{-v}\,dV
\end{equation}
and
\begin{equation}\label{2ndpentropy}
\frac{d^{2}}{dt^{2}}H_{p}(u)=-p\int_{M}\omega^{p-2}\left(|\nabla\nabla v|_{A}^{2}+{\rm Ric}(\nabla v,\nabla v)\right)e^{-v}\,dV,
\end{equation}
\end{lemma}
\proof
Set $u^{p-1}=e^{-v}$, recall the $p$-entropy in \eqref{pentropy},
we have
\begin{equation*}
H_{p}(u)=\int_{M}ve^{-v}\,dV,
\end{equation*}
Using the fact $\nabla u=-\frac{u}{p-1}\nabla v$ in $p$-heat equation \eqref{pheat}, we get
\begin{align*}
-e^{-v}\partial_{t}v&=(p-1)^{p-1}{\rm div}\left[\left|-\frac{u}{p-1}\nabla v\right|^{p-2}\left(-\frac{u}{p-1}\right)\nabla v\right]\\
&=-{\rm div}(e^{-v}|\nabla v|^{p-2}\nabla v)\\
&=-e^{-v}(\Delta_{p}v-|\nabla v|^{p}).
\end{align*}
Let $\omega=|\nabla v|^{2}$, we obtain
\begin{align}\label{pheat11}
\partial_{t}v&=\Delta_{p}v-\omega^{\frac{p}{2}}\notag\\
&=\left(\frac{p}{2}-1\right)\omega^{\frac{p}{2}-2}\langle\nabla \omega,\nabla v\rangle
+\omega^{\frac{p}{2}-1}\Delta v-\omega^{\frac{p}{2}}.
\end{align}
To evaluate the derivatives of $p$-entropy, integration by parts implies
\begin{align*}
\frac{d}{dt}H_{p} &=\frac{d}{dt}\int_{M}ve^{-v}\,dV
=\int_{\mathbb{R}^{n}}(1-v)e^{-v}\partial_{t}v\,dV\\
&=\int_{M}(1-v)e^{-v}(\Delta_{p}v-\omega^{\frac{p}{2}})\,dV\\
&=\int_{M}\omega^{\frac{p}{2}}e^{-v}\,dV,
\end{align*}
and
\begin{align}\label{2ndent}
\frac{d^{2}}{dt^{2}}H_{p}&=\frac{d}{dt}\int_{M}\omega^{\frac{p}{2}}e^{-v}\,dV
\notag\\
&=\frac{p}{2}\int_{M}\omega^{\frac{p}{2}-1}e^{-v}\partial_{t}\omega
\,dV-\int_{V}\omega^{\frac{p}{2}}e^{-v}\partial_{t}v\,dV
\notag\\
&=p\int_{M}\omega^{\frac{p}{2}-1}e^{-v}\nabla v\cdot\nabla
\partial_{t}v\,dV-\int_{M}\omega^{\frac{p}{2}}e^{-v}\partial_{t}v\,dV.
\end{align}

$1)$  Compute the first part of formula \eqref{2ndent}. Since
\begin{align*}
&\nabla v\cdot\nabla\partial_{t}v
=\nabla v\cdot\nabla\left[\left(\frac{p}{2}-1\right)\omega^{\frac{p}{2}-2}\langle\nabla \omega,\nabla v\rangle
+\omega^{\frac{p}{2}-1}\Delta v-\omega^{\frac{p}{2}}\right]
\\
=&\left(\frac{p}{2}-1\right)\left(\frac{p}{2}-2\right)\omega^{\frac{p}{2}-3}\langle\nabla \omega,\nabla
v\rangle^{2}
+\left(\frac{p}{2}-1\right)\omega^{\frac{p}{2}-2}\langle\nabla \omega,\nabla v\rangle\Delta v
\\
+&\omega^{\frac{p}{2}-1}\langle\nabla\Delta v,\nabla
v\rangle-\frac{p}{2}\omega^{\frac{p}{2}-1}\langle\nabla \omega,\nabla v\rangle
+\left(\frac{p}{2}-1\right)\omega^{\frac{p}{2}-2}\bigg(\langle\nabla\nabla v,\nabla
\omega\rangle+\langle\nabla v,\nabla\nabla \omega\rangle\bigg)\cdot\nabla v,
\end{align*}
we rewrite the first part of formula \eqref{2ndent} as the formula \eqref{thefirstpart}
\begin{align}\label{thefirstpart}
&p\left(\frac{p}{2}-1\right)\left(\frac{p}{2}-2\right)\int_{M}\omega^{p-4}e^{-v}\langle\nabla
\omega,\nabla v\rangle^{2}\,dV
\notag
+p\left(\frac{p}{2}-1\right)\int_{M}\omega^{p-3}e^{-v}\langle\nabla \omega,\nabla
v\rangle\Delta v\,dV
\notag\\
& +p\int_{M}\omega^{p-2}e^{-v}\langle\nabla\Delta v,\nabla v\rangle\,dV
-\frac{p^{2}}{2}\int_{M}\omega^{p-2}e^{-v}\langle\nabla \omega,\nabla
v\rangle\,dV
\notag\\
& +p\left(\frac{p}{2}-1\right)\int_{M}\omega^{p-3}e^{-v}\bigg(\langle\nabla\nabla v,\nabla
\omega\rangle+\langle\nabla v,\nabla\nabla \omega\rangle\bigg)\cdot\nabla v\,dV.
\end{align}
For the second term of the formula \eqref{thefirstpart}, integration by parts yields
\begin{align*}
& p\left(\frac{p}{2}-1\right)\int_{M}\omega^{p-3}e^{-v}\langle\nabla\omega,\nabla
v\rangle\Delta vdV
\\
=&-p(p-3)\left(\frac{p}{2}-1\right)\int_{M}\omega^{p-4}e^{-v}\langle\nabla\omega,\nabla
v\rangle^{2}\,dV
+ p\left(\frac{p}{2}-1\right)\int_{M}\omega^{p-2}e^{-v}\langle\nabla\omega,\nabla v\rangle\,dV
\\
- &p\left(\frac{p}{2}-1\right)\int_{M}\omega^{p-3}e^{-v}\bigg(\langle\nabla\nabla\omega,\nabla
v\rangle+\langle\nabla\omega,\nabla\nabla v\rangle\bigg)\cdot\nabla v\,dV.
\end{align*}
For the third term of the formula \eqref{thefirstpart}, one can get
\begin{align*}
& p\int_{M}\omega^{p-2}e^{-v}\langle\nabla\Delta v,\nabla v\rangle\,dV
\\
=&-\frac{p}{2}\int_{M}\nabla(\omega^{p-2}e^{-v})\cdot\nabla\omega\,dV
-p\int_{M}\omega^{p-2}e^{-v}\left(|\nabla\nabla v|^{2}+{\rm Ric}(\nabla v,\nabla v)\right)\,dV
\\
=&-\frac{p(p-2)}{2}\int_{M}\omega^{p-3}e^{-v}|\nabla\omega|^{2}\,dV
+\frac{p}{2}\int_{M}\omega^{p-2}e^{-v}\langle\nabla v,\nabla\omega\rangle\,dV
\\
-&p\int_{M}\omega^{p-2}e^{-v}\left(|\nabla\nabla v|^{2}+{\rm Ric}(\nabla v,\nabla v)\right)\,dV.
\end{align*}
By using of integration by parts and the Bochner formula
\begin{equation*}
\Delta|\nabla v|^{2}-2\nabla v\cdot\nabla\Delta v=2|\nabla\nabla v|^{2}+2{\rm Ric}(\nabla v,\nabla v).
\end{equation*}
Thus, collecting these terms back into formula \eqref{thefirstpart}, we get
\begin{align}\label{part1}
&-\frac{p(p-2)^{2}}{4}
\int_{M}\omega^{p-2}e^{-v}\frac{\langle\nabla\omega,\nabla v\rangle^{2}}{\omega^{2}}\,dV
-\frac{p}{2}\int_{M}\omega^{p-2}e^{-v}\langle\nabla\omega,\nabla v\rangle\,dV
\notag\\
&-\frac{p(p-2)}{2}\int_{M}\omega^{p-2}e^{-v}\frac{|\nabla\omega|^{2}}{\omega}\,dV
-p\int_{M}\omega^{p-2}e^{-v}\left(|\nabla\nabla
v|^{2}+{\rm Ric}(\nabla v,\nabla v)\right)\,dV.
\end{align}

$2)$ Compute the second part of formula \eqref{2ndent}. Recalling the detail expression of $p$-heat equation \eqref{pheat11} and integrating by parts, we have
\begin{align}\label{part2}
&-\int_{M}\omega^{\frac{p}{2}}e^{-v}\partial_{t}v\,dV
\notag\\
=&-\left(\frac{p}{2}-1\right)\int_{\mathbb{R}^{n}}\omega^{p-2}e^{-v}\langle\nabla \omega,\nabla
v\rangle\,dV
-\int_{\mathbb{R}^{n}}\omega^{p-1}e^{-v}\Delta v\,dV+\int_{\mathbb{R}^{n}}\omega^{p}e^{-v}\,dV
\notag\\
=&-\left(\frac{p}{2}-1\right)\int_{M}\omega^{p-2}e^{-v}\langle\nabla \omega,\nabla
v\rangle\,dV+\int_{M}\nabla(\omega^{p-1}e^{-v})\cdot\nabla
v\,dV+\int_{M}\omega^{p}e^{-v}\,dV
\notag\\
=&\frac{p}{2}\int_{M}\omega^{p-2}e^{-v}\langle\nabla \omega,\nabla v\rangle\,dV.
\end{align}
Adding the two parts \eqref{part1} and \eqref{part2} together establishes
\begin{align*}
\frac{d^{2}}{dt^{2}}H_{p}
=&-p\int_{M}\omega^{p-2}e^{-v}\left(\frac{(p-2)^{2}}{4}
\frac{\langle\nabla\omega,\nabla v\rangle^{2}}{\omega^{2}}
+\frac{p-2}{2}\frac{\langle\nabla\omega,\nabla v\rangle^{2}}{\omega^{2}}
+|\nabla\nabla v|^{2}\right)e^{-v}\,dV
\\
&-p\int_{M}\omega^{p-2}e^{-v}{\rm Ric}(\nabla v,\nabla v)\,dV
\\
=&-p\int_{M}\omega^{p-2}\left(|\nabla\nabla
v|_{A}^{2}+{\rm Ric}(\nabla v,\nabla v)\right)e^{-v}\,dV,
\end{align*}
where
\begin{align*}
|\nabla\nabla v|_{A}^{2}
&=a^{ij}a^{kl}v_{ik}v_{jl}
=\left(g^{ij}+(p-2)\frac{v^{i}v^{j}}{w}\right)\left(g^{kl}+(p-2)\frac{v^{k}v^{l}}{w}\right)v_{ik}v_{jl}
\\
&=(p-2)^{2}\frac{(v^{i}v^{k}v_{ik})(v^{j}v^{l}v_{jl})}{\omega^{2}}
+(p-2)\frac{g^{ij}v^{k}v^{l}v_{ik}v_{jl}+g^{kl}v^{i}v^{j}v_{ik}v_{jl}}{\omega}+g^{ij}g^{kl}v_{ik}v_{jl}
\\
&=\frac{(p-2)^{2}}{4}\frac{|\nabla
v\cdot\nabla\omega|^{2}}{\omega^{2}}+\frac{p-2}{2}\frac{|\nabla\omega|^{2}}{\omega}+|\nabla\nabla
v|^{2}.
\end{align*}
\endproof
Now we can start our proof.
\begin{proof}[\bf{Proof of Theorem \ref{Error}}]

The proof need us to evaluate the first and second order derivatives of the $p$-entropy power $N_{p}(u)$.
The first order derivative of the $p$-entropy power is
\begin{equation*}
\frac{d}{dt}N_{p}(u)=\frac{p}{n}N_{p}(u)I_{p}(u),
\end{equation*}
which is calculated by the $p$-DeBruijn's identity defined in $p$-Fisher information \eqref{pfisher}. Define $J_{p}(u)$
\begin{equation}\label{pdeb2}
J_{p}(u)\doteqdot-\frac1p\frac{d}{dt}I_{p}(u)
\end{equation}
as the second order derivative of $p$-entropy, we get
\begin{align*}
\frac{d^{2}}{dt^{2}}N_{p}(u)&=\frac{p}{n}N_{p}(u)\left(\frac{p}{n}
\left(I_{p}(u)\right)^{2}+\frac{d^{2}}{dt^{2}}H_{p}(u)\right)
\\
&=\frac{p^{2}}{n}N_{p}(u)\left(\frac{1}{n}\left(I_{p}(u)\right)^{2}+\frac{1}{p}\frac{d^{2}}{dt^{2}}H_{p}(u)\right)
\\
&=-\frac{p^{2}}{n}N_{p}(u)\left(J_{p}(u)-\frac{1}{n}\left(I_{p}(u)\right)^{2}\right).
\end{align*}
Hence, the concavity condition $\frac{d^{2}}{dt^{2}}N_{p}(u)\leq0$ is
equivalent to
\begin{equation}\label{equiconcavity}
J_{p}(u)\geq\frac{1}{n}\left(I_{p}(u)\right)^{2}.
\end{equation}
Motivated by Villani \cite{Villani}, consider the function $A(\lambda)$
\begin{align}\label{alambda}
A(\lambda)&\doteqdot\int_{M}\left(\left|\omega^{\frac{p}{2}-1}\nabla_{i}\nabla_{j}v+\lambda
a_{ij}\right|_{A}^{2}+\omega^{p-2}{\rm Ric}(\nabla v,\nabla v)\right)e^{-v}\,dV
\notag\\
&=\int_{M}\omega^{p-2}\left(|\nabla\nabla
v|_{A}^{2}+{\rm Ric}(\nabla v,\nabla v)\right)e^{-v}\,dV+n\lambda^{2}+2\lambda\int_{M}e^{-v}\Delta_{p}v\,dV
\notag\\
&=J_{p}(u)+n\lambda^{2}+2\lambda I_{p}(u).
\end{align}
Set $\lambda=-\frac{I_{p}(u)}{n}$, \eqref{alambda} yields the equality
\begin{equation*}
J_p(u)-\frac1n\left(I_p(u)\right)^2=\int_M\left(\left|\omega^{\frac{p}{2}-1}\nabla_{i}\nabla_{j}v
-\frac{I_{p}(u)}{n}a_{ij}\right|_{A}^{2}+\omega^{p-2}{\rm Ric}(\nabla v,\nabla v)\right) e^{-v}\,dV.
\end{equation*}
If the Ricci curvature is nonnegative, we get the inequality \eqref{equiconcavity}.
\end{proof}

\begin{remark}
The concavity inequality \eqref{equiconcavity} can also be proved by Cauchy-Schwarz
inequality. In fact, when Ricci curvature is nonnegative, the $p$-type trace inequality
\begin{equation*}
\omega^{p-2}|\nabla\nabla v|_{A}^{2}\geq\frac{1}{n}\left({\rm tr}_A(\omega^{\frac{p}{2}-1}\nabla_i\nabla_jv)\right)^{2}
=\frac{1}{n}(\Delta_{p}v)^{2}
\end{equation*}
yields
\begin{align*}
J_{p}(u)\ge\frac{1}{n}\int_{M}(\Delta_{p}v)^{2}e^{-v}\,dV
\geq\frac{1}{n}\left(\int_{M}e^{-v}\Delta_{p}v\,dV\right)^{2}
=\frac{1}{n}\left(I_{p}(u)\right)^{2}.
\end{align*}
\end{remark}

\section{Isoperimetric entropy power inequality}

This section is mainly based on the proof of Theorem \ref{isoper}. By observing the first order derivative of the $p$-entropy power, we notice that the concavity property can be interpreted as the decreasing in time property of $t\mapsto\Psi_{p}(u)$ if we introduce a functional $\Psi_{p}(u)$ (defined in \eqref{FI}). Meanwhile, it is a natural idea that one utilizes the monotonicity of a functional of the corresponding solutions to the heat equation equipped with the dilation invariance property to get inequalities in sharp form in Shannon's case. Accordingly, we try to construct inequality \eqref{isoperi} by the above method.

The proof is split into three
parts. In the first part, we verify that $\Psi_{p}(u)$ (defined in \eqref{FI}) is a dilation invariant
function. In the second part, we go further into the monotonicity of the function and
the establishment of the inequality. In the last part, we calculate the value of the
limit.
\begin{proof}[\bf{Proof of Theorem \ref{isoper}}]

Let us first introduce a functional
\begin{equation}\label{FI}
\Psi_{p}(u)\doteqdot N_{p}(u)I_{p}(u),
\end{equation}
Then inequality \eqref{isoperi} in Theorem \ref{isoper} is equivalent to
\begin{equation}\label{Psip}
\Psi_{p}(u)\geq\Psi_{p}\left(\widetilde{G}_{p}(x)\right),
\end{equation}
where $\widetilde{G}_{p}(x)$ defined in formula \eqref{Gp} is the fundamental solution of $p$-heat equation \eqref{pheat}. Therefore, our aim is to obtain the inequality \eqref{Psip}.

Set $\lambda>0$, define a mass-preserving dilation
\begin{equation}\label{dilainvar}
D_{\lambda}u(x)\doteqdot\lambda^{-\frac{n}{p-1}}u(x/\lambda),
\end{equation}
we observe that the functional $\Psi_{p}(u)$ is invariant
\begin{equation}\label{invar}
\Psi_{p}(D_{\lambda}u(x))=\Psi_{p}(u(x)).
\end{equation}
by using the following two identities:
\begin{equation}\label{NPDF}
N_{p}(D_{\lambda}u(x))=\lambda^{p}N_{p}(u(x)),
\end{equation}
\begin{equation}\label{IPDF}
I_{p}(D_{\lambda}u(x))=\lambda^{-p}I_{p}(u(x)).
\end{equation}

We state detailed calculations of the identities below. Recalling the definition of $p$-entropy \eqref{pentropy}, integration by substitution implies
\begin{align*}
H_{p}(D_{\lambda}u(x))&=-\int_{\mathbb{R}^n}(D_{\lambda}u(x))^{p-1}\log (D_{\lambda}u(x))^{p-1}\,dx\\
&=n\log\lambda\int_{\mathbb{R}^{n}}u^{p-1}(x/\lambda)\,d(x/\lambda)
-\int_{\mathbb{R}^{n}}u^{p-1}(x/\lambda)\log
u^{p-1}(x/\lambda)\,d(x/\lambda)
\\
&=n\log\lambda+H_{p}(u(x)),
\end{align*}
i.e.
\begin{equation}\label{HPDF}
H_{p}(D_{\lambda}u(x))=n\log\lambda+H_{p}(u(x)).
\end{equation}
Hence, by the definition of $p$-entropy power \eqref{pSEP}, we have
\begin{equation*}
N_{p}(D_{\lambda}u(x))=\exp\left\{\frac{p}{n}H_{p}(D_{\lambda}u(x))\right\}=\lambda^{p}N_{p}(u(x)).
\end{equation*}
Integration by substitution on the dilation of $u$ of $p$-Fisher information \eqref{pfisher}, we get
\begin{align*}
I_{p}(D_{\lambda}u(x))&=(p-1)^{p}\int_{\mathbb{R}^{n}}\frac{|\nabla
[\lambda^{-\frac{n}{p-1}}u(x/\lambda)]|^{p}}{\lambda^{-\frac{n}{p-1}}u(x/\lambda)}\,dx
\notag\\
&=(p-1)^{p}\lambda^{-p}\int_{\mathbb{R}^{n}}\frac{|\nabla
u(x/\lambda)|^{p}}{u(x/\lambda)}\,d(x/\lambda)
\notag\\
&=\lambda^{-p}I_{p}(u(x)).
\end{align*}
Note that if $\lambda=t^{\frac{1}{p}}$, a direct application of \eqref{NPDF} to \eqref{Gpt} yields \eqref{tGpt}.

The concavity of $p$-entropy power shows that the functional $\Psi_{p}(u)$ is non-increasing on time among solutions to the $p$-heat equation, and it will reach its maximum lower bound as time tends to infinity. Using property as in \eqref{dilainvar}, we scale $u(x,t)$ by this formula
\begin{equation*}
U(x,t)=t^{-\frac{n}{p(p-1)}}u(t^{-\frac{1}{p}}x)
      =D_{t^{\frac{1}{p}}}u.
\end{equation*}
Then $\Psi_{p}(u)=\Psi_{p}(U)$ by \eqref{invar}.
On the other hand, we have $u_{\infty}\doteqdot\lim\limits_{t\rightarrow\infty}U(x,t)=\widetilde{G}_{p}(x).$
Thus, the decreasing in time property of $\Psi_{p}(u)$ implies
\begin{equation}\label{mono}
\Psi_{p}(u)\geq\Psi_{p}(u_{\infty})=\Psi_{p}(\widetilde{G}_{p}(x))=\gamma_{n,p}.
\end{equation}

The last point is the computation of $\gamma_{n,p}$. Let $q=\frac p{p-1}$, using the identities(See \cite{Agueh})
\begin{equation*}
\int_{\mathbb{R}^{n}}e^{-| x|^{q}}\,dx
=\frac{2\pi^{\frac{n}{2}}}{q}\frac{\Gamma(\frac{n}{q})}{\Gamma(\frac{n}{2})}
=\pi^{\frac{n}{2}}\frac{\Gamma(\frac{n}{q}+1)}{\Gamma(\frac{n}{2}+1)},\quad
\int_{\mathbb{R}^{n}}| x|^{q}e^{-| x|^{q}}\,dx =\frac{n}{q}\int_{\mathbb{R}^{n}}e^{-| x|^{q}}\,dx,
\end{equation*}
we get
\begin{equation*}
\int_{\mathbb{R}^{n}}e^{-\frac{p-1}{p^{q}}| x|^{q}}\,dx
=(p^{\frac{1}{p}}q^{\frac{1}{q}})^{n}\pi^{\frac{n}{2}}\frac{\Gamma(\frac{n}{q}+1)}{\Gamma(\frac{n}{2}+1)},
\quad
\frac{p-1}{p^{q}}\int_{\mathbb{R}^{n}}| x|^{q}e^{-\frac{p-1}{p^{q}}| x|^{q}}\,dx
=\frac{n}{q}(p^{\frac{1}{p}}q^{\frac{1}{q}})^{n}
\pi^{\frac{n}{2}}\frac{\Gamma(\frac{n}{q}+1)}{\Gamma(\frac{n}{2}+1)}.
\end{equation*}
Under the above two formulae, compute the $p$-entropy by putting \eqref{Gp} back into \eqref{pentropy}, we have
\begin{align}\label{HpGp}
H_{p}(\widetilde{G}_{p}(x))&=-\int_{\mathbb{R}^{n}}\widetilde{G}_{p}^{p-1}(x)
\log\widetilde{G}_{p}^{p-1}(x)\,dx
\notag\\
&=-C_{p,n}\log C_{p,n}\int_{\mathbb{R}^{n}}e^{-\frac{p-1}{p^{q}}| x|^{q}}\,dx
+\frac{p-1}{p^{q}}C_{p,n}\int_{\mathbb{R}^{n}}| x|^{q}e^{-\frac{p-1}{p^{q}}| x|^{q}}\,dx
\notag\\
&=-\log\left\{(p^{\frac{1}{p}}q^{\frac{1}{q}})^{-n}
\pi^{-\frac{n}{2}}\frac{\Gamma(\frac{n}{2}+1)}{\Gamma(\frac{n}{q}+1)}\right\}
+\frac{n}{q}
\notag\\
&=\log\frac{(p^{\frac{1}{p}}q^{\frac{1}{q}}
e^{\frac{1}{q}})^{n}}{\pi^{-\frac{n}{2}}\frac{\Gamma(\frac{n}{2}+1)}{\Gamma(\frac{n}{q}+1)}}.
\end{align}
Then, the $p$-entropy power of the time-independent function $\widetilde{G}_{p}$ is
\begin{equation}\label{NpGp}
N_{p}(\widetilde{G}_{p}(x))=\exp\left\{\frac{p}{n}H_{p}(\widetilde{G}_{p}(x))\right\}
=p(qe)^{p-1}\pi^{\frac{p}{2}}
\left[\frac{\Gamma(\frac{n}{2}+1)}{\Gamma(\frac{n}{q}+1)}\right]^{-\frac{p}{n}}.
\end{equation}
Note that $\nabla\widetilde{G}_{p}(x)=-\frac{q}{p^{q}}C_{p,n}^{\frac{1}{p-1}}e^{-\frac{| x|^{q}}{p^{q}}}| x|^{q-2}x$, so a direct calculation shows
\begin{equation}\label{IpGp}
I_{p}(\widetilde{G}_{p}(x))=p^{-q}C_{p,n}\int_{\mathbb{R}^{n}}|x|^{q}e^{-\frac{p-1}{p^{q}}|x|^{q}}\,dx
=\frac{n}{p}.
\end{equation}
Therefore, the value of $\gamma_{n,p}=\Psi_{p}(\widetilde{G}_{p}(x))$ conclude in \eqref{gammacon}.
\end{proof}

\section{$L^p$-Euclidean Nash Inequality}
One of the applications based on isoperimetric entropy power inequality \eqref{isoperi} in Theorem \ref{isoper} is to study $L^p$-Nash inequality.
\begin{proof}[\bf{Proof of Theorem \ref{OPNash}}]
Let $u$ be a positive solution to $p$-heat equation on $\mathbb{R}^n$, set $\xi=\| u^{p-1}\|_{L^{1}}\neq1$ and $\zeta^{p-1}(x)=\frac{u^{p-1}(x)}{\xi}$, we observe $\zeta^{p-1}(x)$ is a probability density. If we rewrite the formula $\Psi_{p}(u)\geq\Psi_{p}(\widetilde{G}_{p}(x))$ in \eqref{mono} as
\begin{equation}\label{inequali}
\frac{I_{p}(u)}{I_{p}(\widetilde{G}_{p})}\geq
\exp\left\{-\frac{p}{n}[H_{p}(u)-H_{p}(\widetilde{G}_{p})]\right\},
\end{equation}
then we get
\begin{align*}
I_p(u)&=I_{p}(\xi^{\frac{1}{p-1}}\zeta)=\xi I_{p}(\zeta)
\notag\\
&\geq \xi
I_{p}(\widetilde{G}_{p})\exp\left\{\frac{p}{n}H_{p}(\widetilde{G}_{p})\right\}
\exp\left\{-\frac{p}{n}H_{p}(\zeta)\right\}
\notag\\
&=\xi
I_{p}(\widetilde{G}_{p})\exp\left\{\frac{p}{n}[H_{p}(\widetilde{G}_{p})-\log\xi]\right\}
\exp\left\{-\frac{p}{n}[H_{p}(\zeta)-\log\xi]\right\}
\notag\\
&=\xi
I_{p}(\widetilde{G}_{p})
\exp\left\{\frac{p}{n}\frac{1}{\xi}H_{p}(\xi^{\frac{1}{p-1}}\widetilde{G}_{p})\right\}
\exp\left\{-\frac{p}{n}\frac{1}{\xi}H_{p}(\xi^{\frac{1}{p-1}}\zeta)\right\},
\end{align*}
where we use the identity $ H_{p}(\xi^{\frac{1}{p-1}}\zeta)=\xi H_{p}(\zeta)-\xi\log\xi.$
Applying the formulae \eqref{IpGp} and \eqref{HpGp}, we obtain
\begin{align}\label{infoineq}
I_{p}(u)&\geq\xi\frac{n}{p}\exp\left\{-\frac{p}{n}\frac{1}{\xi}H_{p}(u)\right\}
\exp\left\{\frac{p}{n}H_{p}(\widetilde{G_{p}})\right\}\exp\left\{-\frac{p}{n}\log\xi\right\}
\notag\\
&=\xi^{1-\frac{p}{n}}n(qe)^{p-1}\pi^{\frac{p}{2}}
\left[\frac{\Gamma(\frac{n}{2}+1)}{\Gamma(\frac{n}{q}+1)}\right]^{-\frac{p}{n}}
\exp\left\{-\frac{p}{n}\frac{1}{\xi}H_{p}(u)\right\}.
\end{align}
Let $g(x)$ be a probability density function and $u(x)=g^{q}(x)$, $q=\frac{p}{p-1}$, then
\begin{align*}
 H_{p}(u)&=H_{p}(g^{q})
=-\int_{\mathbb{R}^{n}}g^{(p-1)q}(x)\log
g^{(p-1)q}(x)\,dx\\
&=-q\int_{\mathbb{R}^{n}}\left(g^{p-1}(x)\log g^{p-1}(x)\right)g(x)\,dx.
\end{align*}
Hence, we have
\begin{equation}\label{Jensen}
-H_{p}(g^{q})\geq q\left(\int_{\mathbb{R}^{n}}g^{p}(x)\,dx\right)
\log\left(\int_{\mathbb{R}^{n}}g^{p}(x)\,dx\right)=q\xi\log\xi
\end{equation}
by using Jensen's inequality on convex function $u\mapsto u\log u$ on $\mathbb{R}^{n}_+$.

Combining  inequalities \eqref{infoineq} and \eqref{Jensen} with  the identity
\begin{equation}
I_{p}(u)=I_{p}(g^{q})=p^{p}\int_{\mathbb{R}^{n}}|\nabla g(x)|^{p}\,dx,
\end{equation}
we obtain
\begin{align}\label{PNash}
\left(\int_{\mathbb{R}^{n}}g^{p}\,dx\right)^{1+\frac{q}{n}}
&\leq
\frac{p}{n}\left(\frac{p-1}{e}\right)^{p-1}\pi^{-\frac{p}{2}}
\left[\frac{\Gamma(\frac{n}{2}+1)}{\Gamma(\frac{n}{q}+1)}\right]^{\frac{p}{n}}
\int_{\mathbb{R}^{n}}|\nabla g(x)|^{p}\,dx
\notag\\
&=p^{p}\gamma_{n,p}^{-1}\int_{\mathbb{R}^{n}}|\nabla g(x)|^{p}\,dx,
\end{align}
where $\gamma_{n,p}$ is defined in \eqref{gammacon}.
If $\|g\|_{L^1}\neq1$,  we can obtain the general $L^p$-Euclidean Nash inequality \eqref{PNashineq} by replacing $g$ with $g/\|g\|_{L^1}$ in \eqref{PNash}.
\end{proof}

\section{$L^p$-Euclidean logarithmic Sobolev inequality }
In \cite{Toscani}, G.Toscani showed the sharp logarithmic Sobolev inequality
as a direct result of the concavity of entropy power, moreover, he gave an improvement of the logarithmic Sobolev inequality. In this section, we obtain a new proof of sharp $L^p$-version of logarithmic Sobolev inequality (in \cite{Del Pino} and \cite{Gentil}) applying the concavity of $p$-entropy power, and also get an improvement of $L^p$-logarithmic Sobolev inequality.

Let us first recall the inequality \eqref{inequali}
$$
\frac{I_{p}(u)}{I_{p}(\widetilde{G}_{p})}\geq
\exp\left\{-\frac{p}{n}[H_{p}(u)-H_{p}(\widetilde{G}_{p})]\right\}.
$$
By using of the identities \eqref{HpGp}, \eqref{IpGp} and the facts that $e^{-x}\geq1-x$, we have
\begin{equation}\label{loga}
(p-1)^{p}\int_{\mathbb{R}^{n}}\frac{|\nabla u|^{p}}{u}\,dx
\geq \int_{\mathbb{R}^{n}}u^{p-1}\log u^{p-1}\,dx
+\log\frac{(p^{\frac{1}{p}}q^{\frac{1}{q}}e)^{n}}{\pi^{-\frac{n}{2}}
\frac{\Gamma(\frac{n}{2}+1)}{\Gamma(\frac{n}{q}+1)}},
\end{equation}
the inequality \eqref{loga} is exactly the $L^p$-Euclidean logarithmic Sobolev inequality (in \cite{Del Pino} and \cite{Gentil}). The equivalence of this is given in below.

In fact,
let $u=g^{q},q=\frac{p}{p-1},$ we rewrite the inequality \eqref{loga} as
\begin{equation}\label{logari}
\int_{\mathbb{R}^{n}}g^{p}\log g^{p}\,dx
\leq p^{p}\int_{\mathbb{R}^{n}}|\nabla g|^{p}\,dx
-\log\frac{(p^{\frac{1}{p}}q^{\frac{1}{q}}e)^{n}}{\pi^{-\frac{n}{2}}
\frac{\Gamma(\frac{n}{2}+1)}{\Gamma(\frac{n}{q}+1)}}.
\end{equation}
Changing $g(x)$ into $g_{h}(x)=h^{\frac{n}{p}}g(hx)$, $h>0$, $x\in \mathbb{R}^{n}$, and set $g_{h}(x)$ satisfies
$\int_{\mathbb{R}^{n}}g_{h}^{p}\,dx=1=\int_{\mathbb{R}^{n}}g^{p}\,dx$, we have the following two formulae
\begin{align*}
\int_{\mathbb{R}^{n}}|\nabla g_{h}|^{p}\,dx&=h^{p}\int_{\mathbb{R}^{n}}|\nabla g|^{p}\,dx,\\ \int_{\mathbb{R}^{n}}g_{h}^{p}\log g_{h}^{p}\,dx&=\int_{\mathbb{R}^{n}}g^{p}\log g^{p}\,dx+n\log h.
\end{align*}
Then applying the inequality \eqref{logari} into $g_{h}(x)$ yields
\begin{equation*}
\int_{\mathbb{R}^{n}}g^{p}\log g^{p}\,dx
\leq p^{p}h^{p}\int_{\mathbb{R}^{n}}|\nabla g|^{p}\,dx-n\log h -\log\frac{(p^{\frac{1}{p}}q^{\frac{1}{q}}e)^{n}}{\pi^{-\frac{n}{2}}\frac{\Gamma(\frac{n}{2}+1)}{\Gamma(\frac{n}{q}+1)}}
\end{equation*}
Set $h^{p}=\frac{np^{-(p+1)}}{\int_{\mathbb{R}^{n}}|\nabla g|^{p}\,dx}$, we conclude in
\begin{equation}\label{LpEu}
\int_{\mathbb{R}^{n}}g^{p}\log g^{p}\,dx
\leq\frac{n}{p}\log\left(p^{p}\gamma^{-1}_{n,p}\int_{\mathbb{R}^{n}}|\nabla g|^{p}\,dx\right).
\end{equation}


%

\begin{proof}[\bf Proof of inequality \eqref{GLogSobolev} in Theorem \ref{improve}]
In the above part of this section, we have shown the $L^p$-logarithmic Sobolev inequality \eqref{LpEu} is equivalent to the inequality \eqref{loga}. 
Let us consider the function
\begin{align}\label{HpuHpG}
-H_{p}(u)+H_{p}(\widetilde{G}_{p})
=\int_{\mathbb{R}^{n}}u^{p-1}\log\left(u^{p-1}/\widetilde{G}_{p}^{p-1}\right)\,dx
+\int_{\mathbb{R}^{n}}\left(u^{p-1}-\widetilde{G}_{p}^{p-1}\right)\log\widetilde{G}_{p}^{p-1}\,dx.
\end{align}
If we set $\frac{1}{n}\int_{\mathbb{R}^{n}}|x|^{q}u^{p-1}\,dx\leq p^{\frac{1}{p-1}}$, then the second term in \eqref{HpuHpG} implies
\begin{equation*}
\int_{\mathbb{R}^{n}}\left(u^{p-1}-\widetilde{G}_{p}^{p-1}\right)\log\widetilde{G}_{p}^{p-1}\,dx
=\frac{p-1}{p^{q}}\int_{\mathbb{R}^{n}}|x|^{q}(\widetilde{G}_{p}^{p-1}-u^{p-1})\,dx\geq0,
\end{equation*}
which involves the facts
$\int_{\mathbb{R}^{n}}\widetilde{G}_{p}^{p-1}\,dx=1=\int_{\mathbb{R}^{n}}u^{p-1}\,dx$, and
 $\int_{\mathbb{R}^{n}}| x|^{q}\widetilde{G}_{p}^{p-1}\,dx=np^{\frac{1}{p-1}}$.\\
A Taylor expansion of the first term in \eqref{HpuHpG} yields
\begin{align*}
&\int_{\mathbb{R}^{n}}u^{p-1}\log\left(u^{p-1}/\widetilde{G}_{p}^{p-1}\right)\,dx\\
=&\int_{\mathbb{R}^{n}}\left[\psi\left(u^{p-1}/\widetilde{G}_{p}^{p-1}\right)
-\psi(1)\right]\widetilde{G}_{p}^{p-1}\,dx
\\
\geq&\psi'(1)\int_{\mathbb{R}^{n}}\left(u^{p-1}/\widetilde{G}_{p}^{p-1}-1
\right)\widetilde{G}_{p}^{p-1}\,dx
+\frac{\psi''(1)}{2}\int_{\mathbb{R}^{n}}\left(u^{p-1}/\widetilde{G}_{p}^{p-1}-1
\right)^{2}\widetilde{G}_{p}^{p-1}\,dx
\\
=&\frac{1}{2}\int_{\mathbb{R}^{n}}\left(u^{p-1}/\widetilde{G}_{p}^{p-1}-1
\right)^{2}\widetilde{G}_{p}^{p-1}\,dx,
\end{align*}
where $\psi(x)=x\log x$ is a convex function. For $\theta\in(0,2]$, we obtain
\begin{equation*}
\int_{\mathbb{R}^{n}}\left(u^{p-1}/\widetilde{G}_{p}^{p-1}-1
\right)^{\theta}\widetilde{G}_{p}^{p-1}\,dx
\leq\left(\int_{\mathbb{R}^{n}}\left(u^{p-1}/\widetilde{G}_{p}^{p-1}-1
\right)^{2}\widetilde{G}_{p}^{p-1}\,dx\right)^{\frac{\theta}{2}}
\left(\int_{\mathbb{R}^{n}}\widetilde{G}_{p}^{p-1}\,dx\right)^{1-\frac{\theta}{2}}
\end{equation*}
by H\"older's inequality.
From this, we get the Csiszar-Kullback type inequality
\begin{equation}\label{C-Ktype}
\int_{\mathbb{R}^{n}}u^{p-1}\log\left(u^{p-1}/\widetilde{G}_{p}^{p-1}\right)\,dx
\geq\frac{1}{2}
\left(\int_{\mathbb{R}^{n}}\left(u^{p-1}/\widetilde{G}_{p}^{p-1}-1
\right)^{\theta}\widetilde{G}_{p}^{p-1}\,dx
\right)^{\frac{2}{\theta}}.
\end{equation}
When $\theta=1$, the right hand of inequality \eqref{C-Ktype} takes form in $L^{1}$-norm
\begin{equation*}
\int_{\mathbb{R}^{n}}u^{p-1}\log\left(u^{p-1}/\widetilde{G}_{p}^{p-1}\right)\,dx
\geq\frac{1}{2}|| u^{p-1}-\widetilde{G}_{p}^{p-1}||_{L^{1}}^{2}.
\end{equation*}
When $\theta=1$, $p=2$, \eqref{C-Ktype} is exactly the classic Csiszar-Kullback inequality
used in Shannon's case to reform the logarithmic Sobolev inequality (in \cite{Toscani}).

Going along the way to the improvement, equation \eqref{HpuHpG} indicates
\begin{equation*}
-H_{p}(u)+H_{p}(\widetilde{G}_{p})
\geq\frac{1}{2}\left(\int_{\mathbb{R}^{n}}\left(u^{p-1}/\widetilde{G}_{p}^{p-1}-1\right)
^{\theta}\widetilde{G}_{p}^{p-1}\,dx
\right)^{\frac{2}{\theta}}.
\end{equation*}
Applying the inequality $e^{-x}\geq1-x+\frac{1}{2}x^{2}$ to equation \eqref{inequali},
 we can reform \eqref{loga} as
\begin{align*}
&(p-1)^{p}\int_{\mathbb{R}^{n}}\frac{|\nabla u|^{p}}{u}\,dx
-\left(\int_{\mathbb{R}^{n}}u^{p-1}\log u^{p-1}\,dx
+\log\frac{(p^{\frac{1}{p}}q^{\frac{1}{q}}e)^{n}}{\pi^{-\frac{n}{2}}
\frac{\Gamma(\frac{n}{2}+1)}{\Gamma(\frac{n}{q}+1)}}\right)\\
\geq &\frac{p}{8n}\left(\int_{\mathbb{R}^{n}}\left(u^{p-1}/\widetilde{G}_{p}^{p-1}-1\right)
^{\theta}\widetilde{G}_{p}^{p-1}\,dx
\right)^{\frac{4}{\theta}}.
\end{align*}
Recall the expression of $\gamma_{n,p}$, we know that
\begin{equation*}
\frac{n}{p}\log\left(\frac{pe}{n}\gamma_{n,p}\right)
=\log\frac{(p^{\frac{1}{p}}q^{\frac{1}{q}}e)^{n}}{\pi^{-\frac{n}{2}}
\frac{\Gamma(\frac{n}{2}+1)}{\Gamma(\frac{n}{q}+1)}}.
\end{equation*}
Set $u=g^{q}$, $q=\frac{p}{p-1}$, it implies \eqref{GLogSobolev} (an improvement of $L^p$-logarithmic Sobolev inequality).
\end{proof}

\section*{Acknowledgements}
This work has been partially supported by the National Science Foundation of China, NSFC (No.11701347). The first author would like to thank Professor Xiang-Dong Li for his interest and illuminating discussion.

\end{document}